\documentclass[11pt]{article}
\usepackage{amssymb}
\usepackage{amsmath}
\newtheorem{thm}{Theorem}[section]
\newtheorem{remark}{Remark}[section]

\newtheorem{defn}{Definition}
\newtheorem{proof}{Proof}

\usepackage{graphicx}
\usepackage[all]{xy}
\numberwithin{equation}{section}
\def\R{{\mathbb R}}

\def\P{\operatorname{\bf P}}
\def\E{\operatorname{\bf E}}
\def\U{\operatorname{U}}

\newcommand{\Ha}{\mathbb{H}}
\newcommand{\K}{\mathbb{K}}
\begin{document}
\title{Optimal minimax strategy in a dice game}
\author{F. Crocce\footnote{Facultad de Ciencias, Centro de Matem\'atica. Igu\'a 4225, CP 11400 Montevideo. e-mail: fabian@cmat.edu.uy} \qquad E. Mordecki\footnote{Facultad de Ciencias, Centro de Matem\'atica. Igu\'a 4225, CP 11400 Montevideo. e-mail: mordecki@cmat.edu.uy}}
%\authorone[UdelaR, Montevideo, Uruguay]{Fabi\'an Crocce}
%\address{
%Uruguay}
%. fabian@cmat.edu.uy}
%\authortwo[UdelaR, Montevideo, Uruguay]{Ernesto Mordecki}
%\addresstwo{Facultad de Ciencias, Centro de Matem\'atica. Igu\'a 4225, CP 11400 Montevideo,
%Uruguay. mordecki@cmat.edu.uy}
\date{\today}
\maketitle
\begin{abstract}
Each of two players, by turns, rolls a dice several times accumulating the successive scores until he decides to stop, or he rolls an ace. When stopping, the accumulated turn score is added to the player account and the dice is given to his opponent. If he rolls an ace, the dice is given to the opponent without adding any point. In this paper we formulate this game in the framework of competitive Markov decision processes (also known as \emph{stochastic games}),  show that the game has a value, provide an algorithm to compute the optimal minimax strategy, and present results of this algorithm in three different variants of the game. 
\end{abstract}


%\tableofcontents
\section{Introduction}

Consider a two-players dice game in which players accumulate points by turns with the following rules. The player who reaches a certain fixed number of points is the winner of the game. In his turn each player rolls the dice several times until deciding to stop or rolling an ace. If he decides to stop the accumulated successive scores are added to his account; while if he rolls an ace no additional points are obtained. As a first approach to find optimal strategies for this game Roters \cite{Roters} studied the optimal stopping problem corresponding to the maximisation of the expected score in one turn. The optimal solution is a good way of minimising the number of turns required to reach the objective.

Later, Roters \& Haigh \cite{Haigh} found the strategy that minimises the expected number of turns required to reach the target. This second strategy is better than the one obtained in \cite{Roters} but none of them take into account the consideration of the number of points of the opponent, that is clearly relevant in order to win the game. 

In this paper we formulate this game in the framework of competitive Markov decision processes (also known as \emph{stochastic games}),  show that the game has a value, provide an algorithm to compute the optimal minimax strategy, and present results of this algorithm in three different variants of the game. 

The concept of ``stochastic games'' was introduced in 1953 by Shapley in \cite{Shapley}. In the recent book by Filar and Vrieze \cite{Filar}, the authors provide a general and modern comprehensive approach to this theory departing from the theory of \emph{controlled Markov processes} (that can be considered ``solitaire'' stochastic games) and call the type of games we are interested in \emph{competitive Markov decision processes}, a denomination that we find more accurate than the more usual denomination \emph{stochastic game}.

During the preparation of this paper we found the related article by Neller and Presser \cite{Neller} where the authors, following an heuristic approach, formulate the Bellman equation of the problem (that is a consequence of our results), and compute the optimal strategy of a variant of this game. It must be noted that the theory of Filar and Vrieze \cite{Filar}, that we follow, provide the solution of the problem in the set of all possible strategies, including non-stationary and randomised strategies, i.e. the set of behaviour strategies.

In section 2 we present the theory of \emph{competitive Markov decision processes}, specially in the \emph{transient} case and we conclude the section with the formulation of the theorem we need to solve our dice game. A proof of this theorem can essentially be found in \cite{Filar}. In section 3 we determine the state space of our game, the corresponding action spaces for each player, the payoff function of the game, and the Markov transitions depending on each state of the process and action of the players. In section 4 we present two related games: first, in order to win, the player has to reach the target exactly (if the target is exceeded, he gives the dice to his opponent without changing his accumulated score); in the second variant the players aim to maximise the difference between their scores. In section 5 we present the conclusions.

\section{Competitive Markov decision processes}
A \emph{competitive Markov decision process} (also known as a \emph{stochastic game}) is the mathematical model of a sequential game, in which two players take actions considering the status of a certain Markov processes. Both actions determine an immediate payoff for each player and the probability distribution of the following state of the game. 
Our interest is centred in two-players, finite-state and finite-action, zero-sum games. To define them formally we need the following ingredients:
\begin{itemize}
\item[(S)] States: A finite set $S$ of the possible states of the game.
\item[(A)] Actions: For each state $s \in S$ we consider finite sets $A^s$ and $B^s$ whose elements are the possible actions for the players; at each step both players take his actions simultaneously and independently.
\item[(P)] Payoffs: For each state $s \in S$ a function $r^s:A^s \times B^s \rightarrow \R$ determines the amount that player two has to pay to player one depending on the actions taken by both players.
\item[(TP)] Transition probabilities: For each $(s,a,b)\in S \times A^s \times B^s$, $P_{s,a,b}$ is a distribution probability on $S$, which determines the following state of game.
\end{itemize}
We denote by $S_t$ the state of the game at time $t=0,1,\ldots$; the initial state is fixed $(S_0=s_0)$. At each step $t$ players choose actions $A_t$ and $B_t$ in $A^{S_t}$ $B^{S_t} $ respectively, which determine that player two has to pay to player one an amount of $r^{S_t}(A_t,B_t),$ and the distribution probability of $S_{t+1}$ will be $P_{S_t,A_t,B_t}$.

The random variable $W = \sum_{t=0}^\infty r^{S_t}(A_t,B_t)$ is the total amount that player two pays to player one (could be negative). Note that $W$ depends on the way in which players take their actions. The objective of the game for player one is to maximise the expected value $V$ of the accumulated payoff $W$, while player two has the objective of minimising it. In principle $V$ could be infinite. Often for economic applications the accumulated payoff is $W_{\beta} = \sum_{t=0}^\infty \beta^t r^{S_t}(A_t,B_t)$, called ``the discounted sum'', where $0<\beta<1$ represent the devaluation of the money. This discount factor $\beta$ ensures that $W$ is finite with probability one and the existence of its expected value.

In the case of \emph{transient} stochastic games, the situation considered in this paper, the sum defining $W$ is finite a.s. due to the fact that the process always reaches a final state $s_f \in S$ in a (not necessarily bounded) finite number of steps. In the definition of transient stochastic game additional conditions are required, in order to $V$ to be finite. Before the formal definition of transient stochastic games, the concept of behaviour strategy is introduced. 

\subsection{Strategies}
Consider the set $\K$ defined by
$$\K = \{(s,a,b):s \in S, a\in A^s, b\in B^s \}$$
We define, for each $t=0,1,\ldots$, the sets $\Ha_t$, of admissible histories up to time $t$, by

\begin{equation*}
  \Ha_t=
\left\lbrace
  \begin{array}{ll}
     S &\text{if $t=0$} \\
	 \underbrace{\K \times \ldots \times \K}_{t \text{ times}} \times S &\text{else}
  \end{array}
  \right.
\end{equation*}

\begin{defn}[Behaviour strategy] 
Given a stochastic game, a \emph{behaviour strategy} for player one (two) in the game is a function $\pi$ which associates to each history 
$$h=(s_0,a_0,b_0,\ldots, s) \in \cup_{t=0}^\infty \Ha_t$$
a distribution probability $\pi(\cdot|h)$ in $A^s$ (respectively $\varphi(\cdot|h)$ in $B^s$). In the context of a stochastic game, we denote by $\Pi$ ($\Phi$), the set of all behaviour strategies for player one (two).
\end{defn}

Note that the previous definition is in agreement with the intuitive idea that a player can choose his action based on the history of the game. There are two relevant subclasses of strategies, \emph{pure} and \emph{stationary}, introduced below.

\begin{defn}[Pure strategy] 
A behaviour strategy $\pi$ is said to be \textit{pure}, if for each history $h$ there exists an action $a_h$ such as $\pi(a_h |h)=1$. We could say that a pure strategy chooses the action to be taken in a deterministic way.
\end{defn}

\begin{defn}[Stationary strategy]
A behaviour strategy  $\pi$ is said to be \textit{stationary}, if the probability distribution $\pi(\cdot|h)$ depends only on $s$, the last state of the history. In this case we use the notation $\pi(\cdot|s)$.
\end{defn}

\subsection{Probabilistic framework}

We now construct the probability space in which the optimisation procedure takes place.
%,
Consider the product space 
$$\Omega = (S \times \cup_{s \in S}A^s \times \cup_{s \in S}B^s)^\mathbb{N}$$
equipped with the product $\sigma$-algebra $\mathfrak{F}$, defined as the minimal $\sigma$-algebra containing the cylinder sets of $\Omega$.

Given $\omega = (s_0,a_0,b_0,s_1,a_1,b_1,\ldots) \in \Omega$, a sequence of states and actions in the product space, the coordinate  processes
$$\{S_t\}_{t=0,1,\ldots},\ \{A_t\}_{t=0,1,\ldots},\ \{B_t\}_{t=0,1,\ldots}$$
are defined by
$$S_t(\omega)=s_t,\quad A_t(\omega)=a_t,\quad B_t(\omega)=b_t.$$

In this framework, given $\pi$, $\varphi$ behaviour strategies for players one and two and an initial state $s\in S$, it is possible to introduce a probability $\P_{s,\pi,\varphi}$, such that, for the random vector
$$H_t=(S_0,A_0,B_0,S_1,A_1,B_1,\ldots,S_t),$$
and the finite sequence of states and actions $h_t=(s_0,a_0,b_0,\ldots,s_t),$ the following assertions hold:
\begin{itemize}
\item the game starts in the state $s$, i.e., $\P_{s,\pi,\varphi}(S_0 =s)=1;$
\item with probability 1, $H_t$ take their values in $\Ha_t$;
\item the probability distribution on the actions chosen by players at time $t$ depends on $H_t$, according to
$$\P_{s,\pi,\varphi}(A_t=a_t|H_t=h_t)=\pi(a_t|h_t)$$
$$\P_{s,\pi,\varphi}(B_t=b_t|H_t=h_t)=\varphi(b_t|h_t)$$
$$\P_{s,\pi,\varphi}(A_t=a_t, B_t=b_t|H_t=h_t)=\pi(a_t|h_t)\varphi(b_t|h_t).$$
\item the distribution probability of $S_{t+1}$ depends only on $S_t,A_t,B_t$, being the transition probabilities (TP) of the game
$$\P_{s,\pi,\varphi}(S_{t+1}=s_{t+1}|H_t=h_t, A_t=a_t, B_t=b_t)=P_{s_t,a_t,b_t}(s_{t+1}).$$
\end{itemize}
%Throughout this paper 
We denote by $\E_{s,\pi,\varphi}$ the expected value in the probability space $(\Omega, \mathfrak{F}, \P_{s,\pi,\varphi}).$
\subsection{Transient stochastic games}
\begin{defn}[Transient stochastic game]
\label{deftransient}
A stochastic game is \emph{transient} when there exists a final state $s_f \in S$ such that:
\begin{enumerate}
\item[(1)] $r^{s_f}(a,b)=0,\ \forall a \in A^{s_f},\ \forall b\in B^{s_f};$
\item[(2)] $P_{s_f,a,b}(s_f)=1,\ \forall a \in A^{s_f},\ \forall b\in B^{s_f};$
\item[(3)] for all pair of strategies $(\pi, \varphi)$ of player one and two respectively and for all initial state $s$
$$\sum_{t=0}^\infty \P_{s,\pi,\varphi}(S_t\neq s_f)<\infty.$$
\end{enumerate}
Conditions (1) and (2) ensure that, once the game falls into the final state $s_f$, it never changes the state again and the gain of both players is zero. The third condition ensures that the game finishes with probability one.
\end{defn}

\begin{defn}[Value of a pair of strategies]
Given $\pi \in \Pi,\ \varphi \in \Phi$, strategies for players one and two in a transient stochastic game, the value of the strategies is a function $V_{\pi,\varphi}:S\rightarrow \R$ defined by
$$V_{\pi,\varphi}(s)=\sum_{t=0} ^\infty \E_{s,\pi,\varphi} (r^{S_t}(A_t,B_t)).$$
\end{defn}
\begin{defn}[Optimal strategy]
A behaviour strategy $\pi^*$ for player one in a transient stochastic game is said to be optimal if
$$\inf_{\varphi \in \Phi} V_{\pi^*, \varphi}(s) = \sup_{\pi \in \Pi}\inf_{\varphi \in \Phi}V_{\pi,\varphi}(s), \quad \forall s\in S.$$
Analogously a behaviour strategy $\varphi^*$ for player two, is said to be optimal if
$$\sup_{\pi \in \Pi} V_{\pi, \varphi^*}(s) = \inf_{\varphi \in \Phi}\sup_{\pi \in \Pi}V_{\pi,\varphi}(s), \quad \forall s\in S.$$
\end{defn}
%Rewriting theorem 4.2.6 in \cite{Filar} to this particular case, using our notation we have the following result:
We now formulate the result used to solve our dice game.
\begin{thm}[Value and optimal strategies]
\label{thmoptimal}
Given a transient stochastic game the following identity is fulfilled
\begin{equation}\label{eq:value}
\sup_{\pi} \inf_{\varphi} V_{\pi, \varphi}(s) =  \inf_{\varphi} \sup_{\pi}V_{\pi, \varphi}(s),
\quad
\text{for all $s\in S$.}
\end{equation}
The vector defined in \eqref{eq:value}, denoted by $\left(v(s)\right)_{s\in S}$, is called the \emph{value of the game}. 
This value is the unique joint solution of 
$$
x(s)=\left[r^s(a,b)+\sum_{s'\in S} P_{s,a,b}(s')x(s')\right]^*_{a\in A^s,b\in B^s},
$$
where $[\cdot]^*$ represents the value of the matrix game (in the minimax sense) obtained by considering rows 
$a\in A^s$ and columns $b\in B^s$. 
Moreover, the stationary strategies $\pi$ and $\varphi$ for players one and two, 
%composed by taking at each step, the optimal strategy  
%satisfying 
%that coincides at each step with 
such that $\pi(\cdot|s)$ and $\varphi(\cdot|s)$ are the optimal strategies  
of the matrix game
$$\left[r^s(a,b)+\sum_{s'\in S} P_{s,a,b}(s')v(s')\right]_{a\in A^s,b\in B^s}$$
for all state $s\in S$, are optimal strategies in the transient stochastic game. 
\end{thm}

\begin{proof}
This theorem is essentially a particular case of theorem 4.2.6 in \cite{Filar}. A detailed proof can be found at \cite{Crocce}.
\end{proof}

\begin{remark}
The previous theorem ensures the existence of optimal strategies for both players. Particularly they are in the subclass of stationary strategies. In the proof of this theorem a map $\U$ such that
\begin{equation}
\U v(s)=\left[r^s(a,b)+\sum_{s'\in S} P_{s,a,b}(s')v(s')\right]^*_{a\in A^s,b\in B^s},
\label{eq1}
\end{equation}
which is a $n$-step contraction, is considered. 
Afterwards, we used the map $\U$ to implement a numerical method to find a unique fixed point, that is the value of the game.
\end{remark}

\section{The dice game}
In this section we describe the states, actions, payoffs, and transition probabilities (defined in section 2) corresponding to our dice game, and present the numerical results, showing the optimal strategy for a player.
This strategy, optimal in the class of behaviour strategies, ensures a player to win with probability of at least 1/2 independently of the opponent strategy. The optimal strategy is pure and stationary, and consists in a simple  rule indicating whether \textit{to roll} or \textit{to stop}, depending on the scores of the player and its opponent.
%The information in which the player can base his decision is the course of the game so far.

\subsection{Modelling the dice game}
To solve the dice game (compute optimal strategies), we model it as a transient stochastic game. We have to specify the set of states, possible actions, payoffs and transition probabilities.
\begin{itemize}
\item[(S)] States: During the dice game, there are four aspect varying: the player $j$ who has the dice $(j=1,2)$, the accumulated score $\alpha$ of player one, the accumulated score $\beta$ of player two and the turn score $\tau$ of the player $j$. So, we consider states $(j,\alpha,\beta,\tau)$. We also need to consider two special states: an initial state $s_0$, and a final state $s_f$.

If the score of either of the players is greater or equal than $200$ the game is over, then, it is in the state $s_f$. Because of that, states $(j,\alpha,\beta,\tau)$, only make sense if $\alpha<200$, $\beta<200$. The same happen if $\tau$ is big enough to reach $200$ stopping. So the finite set $S$ of possible states is
$$S=\{s_0,s_f\} \cup S_1 \cup S_2$$
where $S_1$ is the set of states of the player one:
$$S_1=\{(1,\alpha,\beta,\tau): 0 \leq \alpha	\leq 199,\ 0 \leq \beta 	\leq 199,\ 0\leq\tau \leq 205-\alpha\}$$
and $S_2$ is the set of states of the player two:
$$S_2=\{(2,\alpha,\beta,\tau): 0 \leq \alpha	\leq 199,\ 0 \leq \beta 	\leq 199,\ 0\leq\tau \leq 205-\beta\}.$$

\item[(A)] Actions:
We have to specify the set of actions per state for each player. Possible actions in this game are \textit{to roll} and \textit{to stop}. We add an extra action \textit{to wait}, which represents that is not the turn of the player. There are some constraints to ensure the transient condition of the stochastic game: in the states $(1,\alpha,\beta,0)_{\alpha<200}$ does not make sense for player one to take the action \textit{to stop} because there's nothing to loose. If in our model we permit taking the action \textit{to stop} with 0 points in the turn $(\tau = 0)$, is easy to see that there exist a pair of strategies that make the game infinite. The same happen if the action \textit{to roll} is possible when stopping is enough to win. The table \ref{tabla} shows the set of possible actions per state.
\begin{table}
\begin{center}
\caption{Possible actions for each player depending on the state of the game.\label{tabla}}
\begin{tabular}{lcc}
\hline  & player one & player two \\ 
\hline \hline $s_i,\ s_f$ & \textit{to wait} & \textit{to wait} \\ 
\hline $(1,\alpha,\beta,0)$ & \textit{to roll}  & \textit{to wait} \\ 
\hline $(1,\alpha,\beta,\tau)_{0<\tau < 200-\alpha}$ & \textit{to roll}, \textit{to stop}  & \textit{to wait} \\ 
\hline $(1,\alpha,\beta,\tau)_{\tau \geq 200-\alpha}$ & \textit{to stop}  & \textit{to wait} \\ 
\hline $(2,\alpha,\beta,0)$ & \textit{to wait}  & \textit{to roll} \\ 
\hline $(2,\alpha,\beta,\tau)_{0<\tau < 200-\beta}$ & \textit{to wait} & \textit{to roll}, \textit{to stop} \\ 
\hline $(2,\alpha,\beta,\tau)_{\tau \geq 200-\beta}$ & \textit{to wait} & \textit{to stop} \\
\hline 
\end{tabular}
\end{center}
\end{table}

\item[(P)] Payoffs:
Because we want to maximise the probability of winning, we define the payoff function 
in such a way that
%so that to 
maximising the probability of winning is equivalent to maximising the expected value $V$ of the payoffs accumulated along the game. 
The model of transient stochastic game allow us to define a payoff for each pair $(state,action)$ but in this case is enough to define the payoff depending only on the state as follows:
\begin{equation*}
  r^s=
\left\lbrace
  \begin{array}{ll}
     1 &\text{if $s=(1,\alpha,\beta,\tau)$  with $\alpha+\tau\geq 200$} \\
	 0 &\text{else}
  \end{array}
  \right.
\end{equation*}
\item[(TP)] Transition probabilities:
To represent graphically the transition probabilities we use the following representation for the states:
\vspace{-0.9em}
$$(1,\alpha,\beta,\tau)=\xymatrix{*++[o][F]{\alpha^\tau_\beta}} \qquad (2,\alpha,\beta,\tau) = \xymatrix{*+[F]{\beta^\tau_\alpha}}.$$
The dynamic of the game and the semantic of the states determine transition probabilities between states. 
In the figure \ref{figtrans}, the probability transitions from a state $(1,\alpha,\beta,\tau)$ with $\alpha+\tau<200$, depending on the player decision are presented.\\

\begin{figure}[hbt]
\begin{center}
$
\xymatrix@=0.5cm{
	&\text{\textit{to roll}}&&\\
	& *++[o][F]{\alpha^\tau_\beta} \ar[drr]^{\frac{1}{6}} \ar[dl]_{\frac{1}{6}} \ar[ddrr]^{\frac{1}{6}} \ar[ddl]^{\frac{1}{6}}    \ar[dd]^{\frac{1}{6}} \ar[ddr]^{\frac{1}{6}}& &\\
	*++[F]{\beta^0_\alpha}& & &  *++[o][F]{\alpha^{\tau+2}_\beta} \\
  	*++[o][F]{\alpha^{\tau+3}_\beta} & *++[o][F]{\alpha^{\tau+4}_\beta}	& *++[o][F]{\alpha^{t+5}_\beta} & *++[o][F]{\alpha^{\tau+6}_\beta}
}
$\qquad
$
\xymatrix@=0.5cm{
	\text{\textit{to stop}}\\
	*++[o][F]{\alpha^\tau_\beta} \ar[dd]^1\\
	\\
	*+[F]{\beta^0_{\alpha+\tau}}
}
$
\end{center}
\caption{Transition probabilities from a state $(1,\alpha,\beta,\tau)$ with $\alpha+\tau<200$, depending on the player decision.\label{figtrans}}
\end{figure}
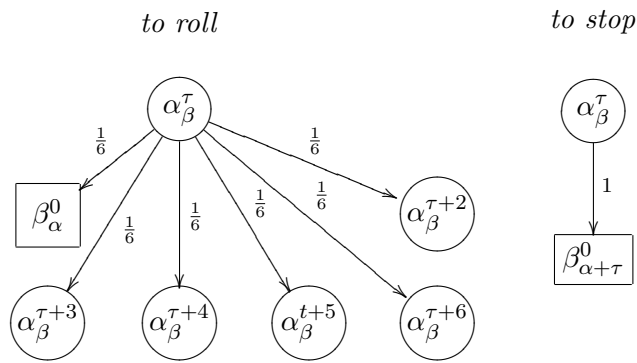

Figure \ref{figtrans} shows that, when the decision is \textit{to roll}, the distribution probability on the states is associated with the results of rolling a dice; particularly the probability of loosing the turn is 1/6. In the winning states of player one (i.e. $(1,\alpha,\beta,\tau)$ with $\alpha+\tau \geq 200$) the transition is, with probability one, to the final state ($s_f$). Transitions for player two are completely symmetric. As is showed in figure \ref{figtransfinal} in the special states $s_i$ and $s_f$ transitions do not depend on the actions taken by players, indeed they do not have options.

\begin{figure}[hbt]
\begin{center}
$
\xymatrix@=0.5cm{
	& *++[o][F]{s_i} \ar[dl]^{\frac{1}{2}} \ar[dr]_{\frac{1}{2}} &\\
	*+[F]{0^0_0} & &  *++[o][F]{0^0_0}
}
$\qquad
$
\xymatrix@=0.5cm{
	*+[o][F]{s_f} \ar@(d,r)_{1}
}
$
\end{center}
\caption{Transition probabilities from the initial and final states.\label{figtransfinal}}
\end{figure}
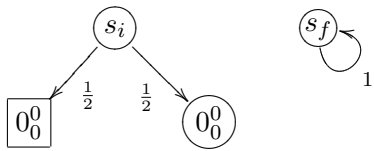

\end{itemize}

Now we verify that the stochastic game defined above is transient. 
We then prove that $s_f$ satisfies conditions (1), (2) and (3) in definition \ref{deftransient}. 
Conditions (1) and (2) are trivially fulfilled, it remains to verify that
$$
\sum_{t=0}^\infty \P_{s,\pi,\varphi}(S_t\neq s_f)<\infty,
$$
for all initial state $s$, and every pair of strategies $\pi$, $\varphi$. 
%The way in which the possible actions in a state were defined does not allow avoid the finish of the game indefinitely, since 

Due to the fact that 
at the beginning of a turn the only option is \textit{to roll}, 
and that in a state in which the accumulated score is enough to win the player has \textit{to stop},
the game can not continue indefinitely.

%In virtue of that, if, f
For example, if a 6 is rolled consecutively 70 times it is impossible to avoid reaching the final state. 
Denoting by $\gamma$ the probability of rolling 70 times a 6, i.e. $\gamma = 1/6^{70}>0$, 
it is easy to see that $\gamma$ is a lower bound for the probability of $S_t=s_f$ for $t\geq 70$. 
By a similar argument, it follows that, for $n=0,1,\ldots$
$$
\P(S_t\neq s_f)<(1-\gamma)^ n\quad\mbox{ if }\quad 70\,n\leq t < 70\,(n+1).
$$
Then
$$\sum_{t=0}^\infty \P_{s,\pi,\varphi}(S_t\neq s_f)< 70 \sum_{n=0}^\infty (1-\gamma)^n < \infty$$
and our model is transient.
\subsection{Numerical results}
\begin{figure}[h]
$
\begin{array}{cc}
\beta=0&\beta=150\\
\includegraphics[scale=0.32]{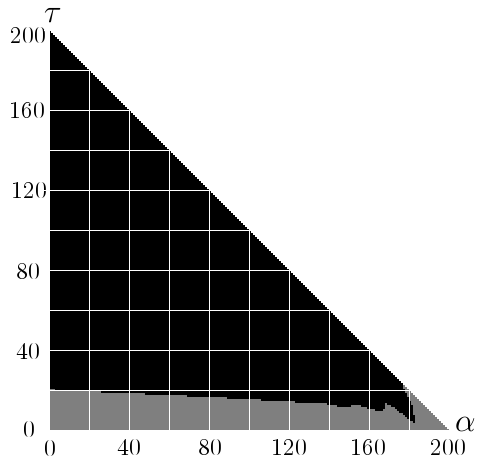}&\includegraphics[scale=0.32]{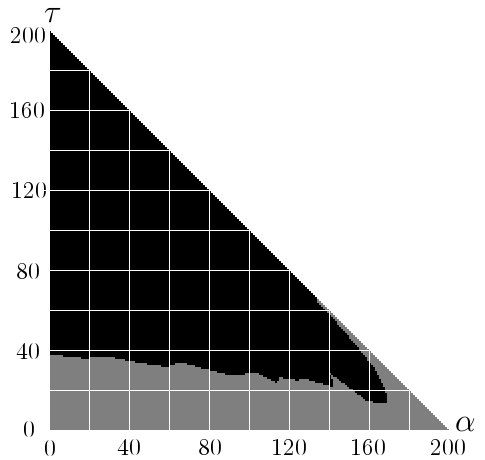}\\ \\

\beta=180&\beta=185\\
\includegraphics[scale=0.32]{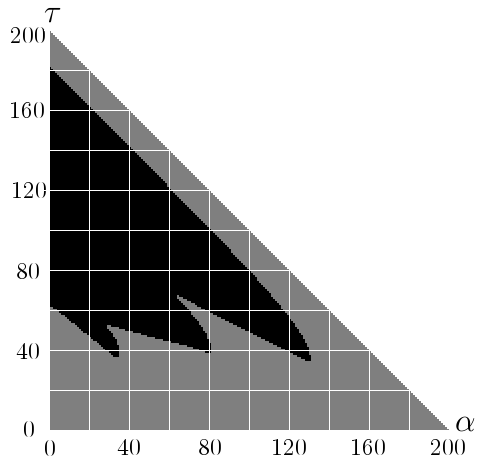}&\includegraphics[scale=0.32]{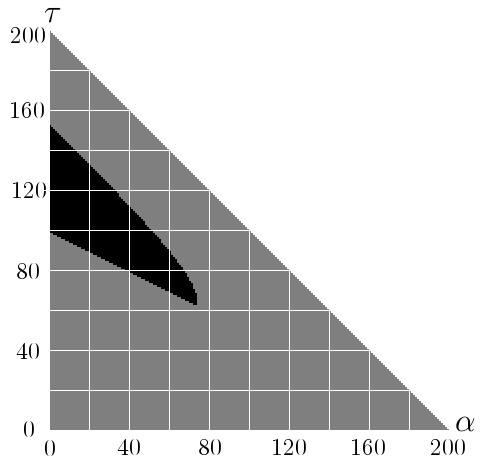}
\end{array}
$
\caption{\label{optima} Part of the optimal strategy for player one (for player two is symmetric). It includes states with opponent score $\beta=0,\ 150,\ 180\ \& \ 185$. In the gray zone the optimal action is \textit{to roll} and in the black zone is \textit{to stop}.}
\end{figure}

In this section the results of the theorem \ref{thmoptimal} are applied  to the particular case of the transient stochastic game defined above. Rewriting the definition of application $\U$, defined in equation \eqref{eq1}, we obtain:

\begin{normalsize}\begin{equation*}
 \U v(s)=
\left\lbrace
  \begin{array}{ll}
     \frac{1}{2}v(1,0,0,0)+\frac{1}{2}v(2,0,0,0)& \text{if $s=s_i$}\\
     v(1,\alpha,\beta,0)_{roll} & \text{if $s=(1,\alpha,\beta,0)$}\\
     \max \left\{ v(1,\alpha,\beta,\tau)_{stop},v(1,\alpha,\beta,\tau)_{roll} \right\}& \text{if $s=(1,\alpha,\beta,\tau):\ \alpha+\tau<200$}\\
	 1 & \text{if $s=(1,\alpha,\beta,\tau):\ \alpha+\tau\geq 200$}  \\
     v(2,\alpha,\beta,0)_{roll} & \text{if $s=(2,\alpha,\beta,0)$}\\
     \min \left\{ v(2,\alpha,\beta,\tau)_{stop},v(2,\alpha,\beta,\tau)_{roll} \right\}& \text{if $s=(2,\alpha,\beta,\tau):\ \beta+\tau<200$}\\
	 0 &\text{in other case}
  \end{array}
  \right.
\end{equation*}\end{normalsize}
where
$$
\begin{array}{lll}
v(1,\alpha,\beta,\tau)_{roll} &=& \frac{1}{6}v(2,\alpha,\beta,0) +  \frac{1}{6}\sum_{k=2}^6 v(1,\alpha,\beta,\tau+k)\\
v(1,\alpha,\beta,\tau)_{stop} &=& v(2,\alpha+\tau,\beta,0)\\
v(2,\alpha,\beta,\tau)_{roll} &=& \frac{1}{6}v(1,\alpha,\beta,0) +  \frac{1}{6}\sum_{k=2}^6 v(2,\alpha,\beta,\tau+k)\\
v(2,\alpha,\beta,\tau)_{stop} &=& v(1,\alpha,\beta+\tau,0).\\
\end{array}
$$
Note that in the equations above we have replaced the value of the matrix games in equation \eqref{eq1} by a maximum, 
in the states in which player one has to take the decision, 
and by a minimum when is player two the one that has to do it. 
In the states in which both players have only one choice the value of the matrix game is the only entry of the matrix.
Since there are no states in which both players have to decide simultaneously, 
the stationary strategy that emerges from the theorem turns out to be pure,
i.e. each player takes an action with probability 1. 
To determine the complete solution is necessary to specify which action should be taken in about $4\,000\,000$ states. 
In figure \ref{optima} the optimal strategy for player one for some states is shown. 
The complete solution can be found at

\verb|www.cmat.edu.uy/cmat/docentes/fabian/documentos/optimalstrategy.pdf|.

\textbf{Some observations about the solution:}
\begin{itemize}
\item At the beginning of the game, when both players have low scores, 
we see that 
the optimal action is \textit{to roll} when $\tau < 20$ and \textit{to stop} in the other case, 
following the strategy found by Roters \cite{Roters} maximising the expected value of a turn score. 

The heuristic interpretation of this fact is: when far away from the target
it is optimal to approach it in steps as large as possible.
\item When the opponent score $\beta$ becomes larger the optimal strategy becomes riskier. This can be explained because there are less turns to reach the target.
\item For opponent scores greater or equal than 187 $(\beta \geq 187)$, the graphic becomes absolutely gray. In other words, when the opponent is  close to win, giving him the dice is a bad idea.
\item To compare the optimal strategy with the one found by Haigh \& Roters \cite{Haigh}, 
we simulate $10\,000\,000$ games. 
Our simulation showed that in $52\%$ of the games, the winner was the player with the optimal strategy.
\end{itemize}

\section{Two related games}
\subsection{Reaching exactly the target}
\label{subsecrebote}
\begin{figure}[htb]
\begin{center}
$
\xymatrix@=0.5cm{
	&\text{\textit{to roll}}&&\\
	& *+++[o][F]{180^{16}_\beta} \ar[ddl]_{\frac{3}{6}} \ar[dd]_{\frac{1}{6}} \ar[ddr]_{\frac{1}{6}} \ar[ddrr]^{\frac{1}{6}}& &\\
	 \\
  	*++[F]{\beta^0_{180}}& *+++[o][F]{180^{18}_\beta} & *+++[o][F]{180^{19}_\beta}	& *+++[o][F]{180^{20}_\beta} & \\ \\
}
$
\end{center}
\caption{Example of transition probability in the variant of the game presented in section \ref{subsecrebote}.\label{figrebote}}
\end{figure}
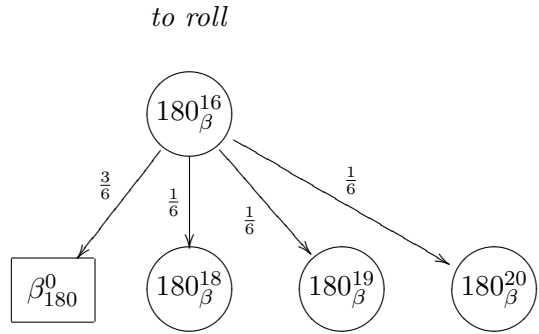

\begin{figure}[h]
$
\begin{array}{cc}
\beta=0&\beta=150\\
\includegraphics[scale=0.32]{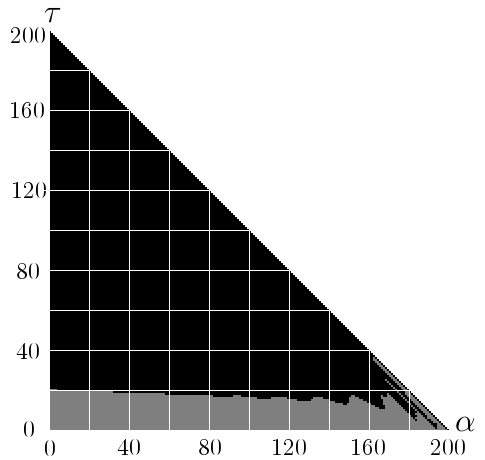}&\includegraphics[scale=0.32]{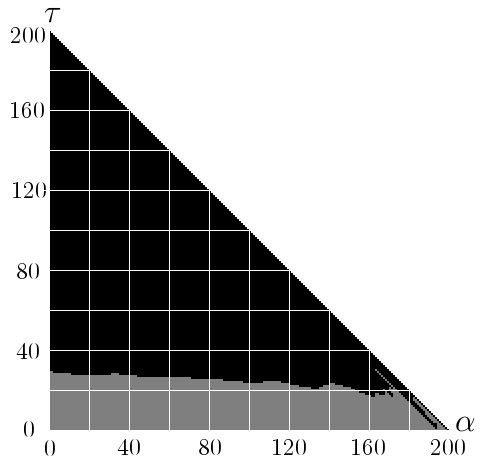}\\ \\

\beta=180&\beta=198\\
\includegraphics[scale=0.32]{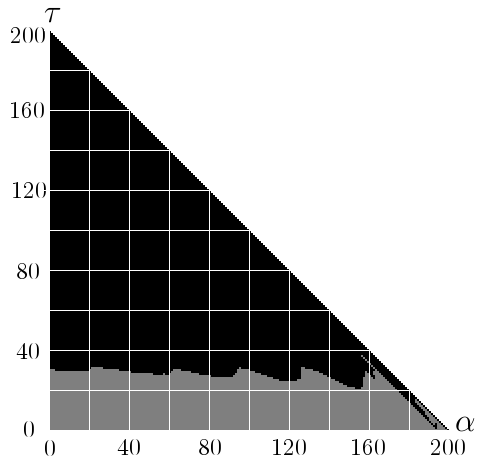}&\includegraphics[scale=0.32]{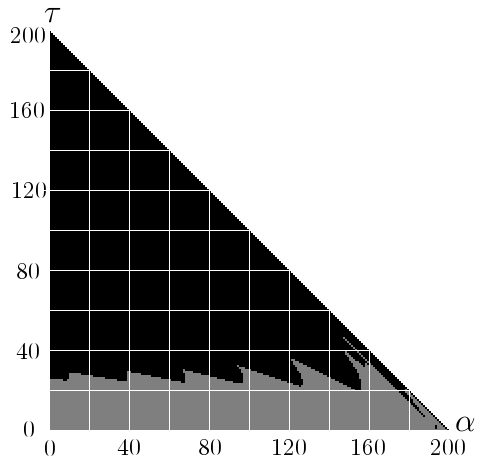}
\end{array}
$
\caption{\label{rebote} Part of the optimal strategy for player one (for player two is symmetric) for the game presented in subsection \ref{subsecrebote}. It includes states with opponent score $\beta=0,\ 150,\ 180\ \& \ 198$. In the gray zone the optimal action is \textit{to roll} and in the black zone is \textit{ to stop}.}
\end{figure}

It is interesting to explore how the optimal strategy changes when the game is modified. 
In this section we consider the same dice game, 
with the only difference being that the condition to win is to reach exactly 200 points. 
If the sum of accumulated and turn score is greater than 200 the turn finishes without capitalising any point. 
%s earned in the current turn.
The formulation of the game is quite similar, 
the difference appears when the accumulated score plus the turn score is greater that 194, situation in which
%is in the states in which the accumulated score of the current player added to the turn score is so close to 200 that with only 
one roll of the dice can exceed the target. 
As an example of the mentioned difference, in figure \ref{figrebote} we show the transition probabilities  
when the accumulated score is 180 and the turn score is 16 ($180+16>194$).
Note that the probability of losing the turn is the probability of rolling a 1,5,6. In figure \ref{rebote} part of the optimal strategy for this variant of the game is shown. The complete optimal strategy is available in  \verb|www.cmat.edu.uy/cmat/docentes/fabian/documentos/optimalexactly.pdf|.

\textbf{Some remarks about the solution:}
\begin{itemize}
\item As in the classical game, when the target is far, the strategy is similar to ``stop if $\tau >20$''.
\item Unlike the classical game, the optimal strategy in this case never becomes so risky. 
This is easy to understand because the probability of winning in any turn is less or equal than $1/6$, 
even in the case of being very close to the target.
%
%independently of the score of the opponent.
%
%despite how close is the opponent to the target score, 
%the probability of winning is not so high, because he need an specific outcome of the dice.
\item When $\alpha + \tau = 194$ there is a ``roll zone'' larger than usual, because 194 is the largest 
score in which there is no risk of losing in one roll but it is possible to win rolling a 6.
\end{itemize}

\subsection{Maximising the expected difference}
\label{subsecmaxdif}
In the second variant of the game,
% we consider, 
%consists 
%in considering that the objective, 
%instead of trying to win, 
the winner, when reaching the target, 
obtains from the loser the difference between the target and the loser's score.

\begin{figure}[htb]
$
\begin{array}{cc}
\beta=0&\beta=150\\
\includegraphics[scale=0.32]{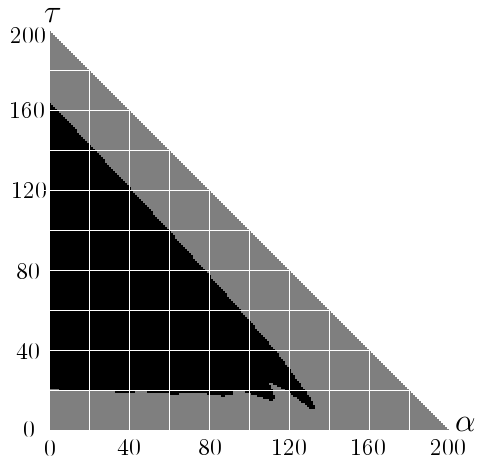}&\includegraphics[scale=0.32]{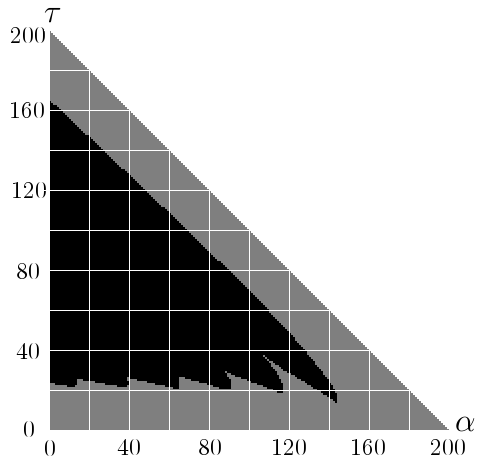}\\ \\
\beta=170&\beta=180\\
\includegraphics[scale=0.32]{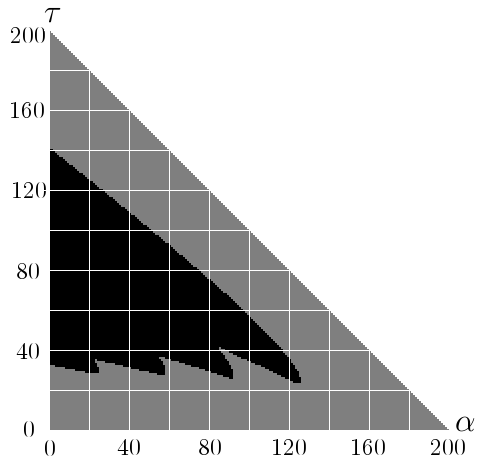}&\includegraphics[scale=0.32]{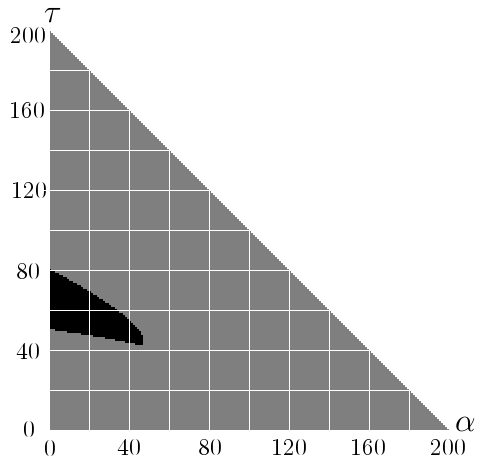}
\end{array}
$

\caption{\label{maxdif}  Part of the optimal strategy for player one (for player two is symmetric) for the variant presented in subsection \ref{subsecmaxdif}. It includes states with opponent score $\beta=0,\ 150,\ 170\ \& \ 180$. In the gray zone the optimal action is \textit{to roll} and in the black zone is \textit{to stop}.}
\end{figure}

%to maximise the difference in scores, 
%assuming that the loser has to pay this difference to the winner.
%maximising the probability of winning, 
%the player has to maximise the expected value of the difference in scores. 
%This variant 
%For example, if when player one wins, player two has 160 points, then he has to pay player one \$ 40. 
Again, the model of the game is very similar to the classical model, 
the only difference is the payoff function:
\begin{equation*}
  r^s=
\left\lbrace
  \begin{array}{ll}
     200-\beta, &\text{if $s=(1,\alpha,\beta,\tau)$  with $\alpha+\tau\geq 200$,} \\
     \alpha-200, &\text{if $s=(2,\alpha,\beta,\tau)$  with $\beta+\tau\geq 200$,} \\
	 0, &\text{else.}
  \end{array}
  \right.
\end{equation*}
Figure \ref{maxdif} shows the optimal strategy for some opponent scores. 
The complete optimal strategy can be found at \\ \verb|www.cmat.edu.uy/cmat/docentes/fabian/documentos/optimalmaxdif.pdf|.\\
%It is more difficult to understand why the optimal strategy is such. 
The main difference with the optimal strategy in the classical case, 
is that when one player is close to win
(taking into account his current turn score), 
he takes the risk of rolling, this feature being observed for any score of the opponent.
%he $\tau$ is large enough to be close to win; unlike the original game, in this cases, the optimal action is \textit{to roll}.

\section{Conclusion}
In this paper we model a dice game in the framework of Markov competitive decision processes (also known as \emph{stochastic games}) 
in order to obtain optimal strategies for a player. Our main results are the proof of the existence of a value and an optimal minimax 
strategy for the game, and the proposal of an algorithm to find this strategy. 
We base our results on the theory of transient stochastic games exposed by Filar and Vrieze in \cite{Filar}.

Previous mathematical treatments of this problem include the solution of the optimal stopping problem for a player that wants to maximise the expected number of points in a single turn (see Roters \cite{Roters}) and the minimisation of the expected number of turns required to reach a target (see Haigh and Roters and \cite{Haigh}). Another previous contribution was done by 
Neller and Presser \cite{Neller}, who found the optimal strategy in the set of stationary pure strategies, departing from a Bellman equation.

We also provide an algorithm to compute explicitly this optimal strategy (that coincides with the optimal strategy in the larger class of behaviour strategies) and show how this algorithm works in three different variants of the game.
\section*{Acknowledgements}
This work was partially supported by Antel-Fundaciba agreement ``An\'alisis de Algoritmos de Codificaci\'on y Cifrado''
\bibliographystyle{apt}

%\bibliography{references_temp}

\begin{thebibliography}{99}
\footnotesize

\bibitem{Crocce}
{\sc Crocce, F.} (2009).
\newblock {\em {Juegos estoc\'asticos transitorios y aplicaciones.}}
\newblock {Master thesis, PEDECIBA, Montevideo, Uruguay.}

\bibitem{Filar}%book
{\sc Filar, J. \& Vrieze, K.} (1997).
\newblock {\em {Competitive Markov Decision Processes.}}
\newblock {Springer-Verlang, New York.}

\bibitem {Haigh}
{\sc Haigh, J. \& Roters, M. }(2000).
\newblock{Optimal Strategy in a Dice Game.}
\newblock {\em Journal of Applied Probability} {\bf 37,} 1110--1116.

\bibitem {Neller}
{\sc Neller, T. \& Presser, C. }(2004).
\newblock{Optimal Play of the Dice Game Pig.}
\newblock {\em The UMAP Journal} {\bf 25.1,} 25--47.

\bibitem {Roters} 
{\sc Roters, M. }(1998). 
\newblock{Optimal Stopping in a Dice Game.}
\newblock {\em Journal of Applied Probability} {\bf 35,} 229--235.

\bibitem {Shapley}
{\sc Shapley, L.S. }(1953).
\newblock{Stochastic games.}
\newblock{\em Proc. of the Nat. Acad. Sciences} {\bf 39,} 1095-1100.
\end{thebibliography}

\end{document}